\documentclass[a4paper,12pt]{article}
\usepackage{amssymb}

\newcommand{\R}{\mathbb R}

\newtheorem{lemma}{Lemma}

\newtheorem{theorem}{Theorem}

\begin{document}

\author{Octavian G. Mustafa\\
\small{Department of Mathematics {\&} Computer Science,}\\
\small{\c{C}ankaya University,}
\small{\"{O}gretmenler Cad. 14 06530 Balgat, Ankara, Turkey}\\
\small{e-mail: octawian@yahoo.com}\\
\small{and}\\
Yong Zhou\\
\small{Department of Mathematics}\\
\small{Zhejiang Normal University, Jinhua, Zhejiang 321004, China}\\
\small{e-mail: yzhoumath@zjnu.cn}}
\title{On the positive solutions to some quasilinear elliptic partial differential equations}
\date{}
\maketitle

\noindent\textbf{Abstract.} We establish that the elliptic equation $\Delta u+f(x,u)+g(\left\vert x\right\vert )x\cdot \nabla u=0$, where
$x\in\mathbb{R}^{n}$, $n\geq3$, and $\vert x\vert>R>0$, has a positive solution which decays to $0$ as $\vert x\vert\rightarrow +\infty$ under
mild restrictions on the functions $f,g$. The main theorem extends and complements the conclusions of the recent paper [M. Ehrnstr\"{o}m, O.G.
Mustafa, On positive solutions of a class of nonlinear elliptic equations, Nonlinear Anal. TMA 67 (2007), 1147--1154]. Its proof relies on a
general result about the long-time behavior of the logarithmic derivatives of solutions for a class of nonlinear ordinary differential equations
and on the comparison method.



\section{Introduction}
This note, motivated by the recent papers \cite{DMM,EM}, is concerned with the existence of a positive solution to the boundary value problem
\begin{eqnarray}
\left\{
\begin{array}{ll}
u^{\prime\prime}+F(t,u)=0,\qquad t\geq t_0>0,\\
q_{-}(t)\cdot\frac{u(t)}{t}\leq u^{\prime}(t)\leq q_{+}(t)\cdot\frac{u(t)}{t},\qquad t\geq t_0,\\
u(t)=o(t)\qquad\mbox{as }t\rightarrow+\infty,
\end{array}
\label{main_bvp} \right.
\end{eqnarray}
for a certain class of continuous functions $F:[t_0,+\infty)\times\mathbb{R}\rightarrow[0,+\infty)$. Here, the functions
$q_{\pm}:[t_0,+\infty)\rightarrow[0,1]$ are assumed continuous.

The interest in studying the $q_{\pm}$ -- problem (\ref{main_bvp}) comes from an investigation of the existence and decay rates of the positive,
vanishing at $+\infty$, solutions to the quasilinear elliptic equation of second order
\begin{eqnarray}
\Delta u+f(x,u)+g(\left\vert x\right\vert )x\cdot \nabla u=0,\qquad x\in G_R,\label{main_PDE}
\end{eqnarray}
where $G_R=\{x\in {\mathbb R}^n:\left\vert x\right\vert
>R\}$ and $n\geq 3$. For an account of recent literature on this topic, we refer to the studies \cite{Constantin1996,Constantin1997,Deng,E,Hesaaraki,Orpel,Wahlen,Yin}.

Following \cite{Constantin1996,NS}, we consider that the functions $f:\overline{G}_{R}\times {\mathbb R}\rightarrow {\mathbb R}$ and
$g:[R,+\infty )\rightarrow {\mathbb R}$ are locally H\"{o}lder continuous. Moreover,
\begin{eqnarray}
0\leq f(x,U)\leq m(\vert x\vert,U),\qquad x\in G_R, \thinspace U\in[0,\varsigma],  \label{core}
\end{eqnarray}
for some $\varsigma>0$ and the continuous application $m:[R,+\infty)\times[0,\varsigma]\rightarrow[0,+\infty)$. The regularity assumptions upon
$f,g$ are sufficient for applying the comparison method \cite{GilbargTrudinger} to the analysis of (\ref{main_PDE}). In fact, given $u(t)$ a
positive solution of (\ref{main_bvp}), the function
\begin{eqnarray}
U(x)=U(\vert x\vert)=\frac{u(t)}{t},\qquad\mbox{where }\vert x\vert=\theta(t)=\left(\frac t{n-2}\right)^{\frac 1{n-2}} \label{change}
\end{eqnarray}
and $t\geq t_0=(n-2)R^{n-2}$, will be a super-solution to (\ref{main_PDE}) satisfying the additional restriction
\begin{eqnarray}
x\cdot\nabla U(x)\leq0,\qquad x\in G_{R}.\label{add_cond_nabla_u}
\end{eqnarray}

It has been noticed in \cite{E} that, when $g$ takes only nonnegative values, the additional requirement (\ref{add_cond_nabla_u}) for the
solution $U$ of the elliptic partial differential equation
\begin{eqnarray}
\Delta U+m(\vert x\vert,U)=0,\qquad x\in G_R,\label{elliptic_m}
\end{eqnarray}
allows for a complete removal of the integral conditions regarding $g$ from the hypotheses of various theorems in the recent literature. Further
developments of this observation are given in \cite{DMM,EM}.

Condition (\ref{add_cond_nabla_u}), translated into the language of ordinary differential equations, reads as
\begin{eqnarray}
u^{\prime}(t)-\frac{u(t)}{t}\leq0,\qquad t\geq t_0.\label{cond_Ehrnstrom}
\end{eqnarray}
It is obvious now that the $q_{\pm}$ -- problem (\ref{main_bvp}) constitutes an improvement of (\ref{cond_Ehrnstrom}).

The use of (\ref{main_bvp}), in the particular case of $q_{-}(t)=0$, $q_{+}(t)=(\ln t)^{-1}$ throughout $[t_0,+\infty)$, has been observed in
\cite{DMM}. To give it a self-contained presentation, let us recall first the essence of the \textit{reduction technique} of \cite{E}. The
comparison ordinary differential equation in the study of super-solutions of (\ref{main_PDE}) (that is, a rewriting of $\Delta U+m(\vert
x\vert,U)+g(\vert x\vert)x\cdot\nabla U=0$, $\vert x\vert>R$, which takes into account (\ref{change})) being displayed as
\begin{eqnarray}
u^{\prime\prime}+H(t,u)+h(t)\left(u^{\prime}-\frac{u}{t}\right)=0,\qquad t\geq t_0,\label{gen_comp}
\end{eqnarray}
where
\begin{eqnarray}
H(t,u)=\frac{1}{n-2}\theta(t)\theta^{\prime}(t)m\left(\theta(t),\frac{u}{t}\right),\quad
h(t)=\theta(t)\theta^{\prime}(t)g(\theta(t)),\label{def_H}
\end{eqnarray}
the method in \cite{E} consists of removing the quantity "$h(t)\left(u^{\prime}-\frac{u}{t}\right)$" from (\ref{gen_comp}) whenever $g$ is
nonnegative-valued and the problem (\ref{main_bvp}), with $q_{-}(t)\equiv0$, $q_{+}(t)\equiv1$, has a positive solution. The approach allows for
total freedom of $g$, however, it keeps intact all the restrictions concerned with $m(\vert x\vert,U)$. Instead of this, we can use the next
modification of (\ref{gen_comp}), namely the quasilinear ordinary differential equation
\begin{eqnarray}
u^{\prime\prime}+\left\vert \frac{H(t,u)}{u}-[1-q_{+}(t)]\cdot\frac{h(t)}{t}\right\vert u=0,\qquad t\geq t_0.\label{new_comp}
\end{eqnarray}
In the new setting (\ref{new_comp}), the functional quantity $\frac{H(t,u)}{u}$ is controlled only partially by the hypotheses of various
comparison-type results.

We also notice that, since
\begin{eqnarray*}
\frac{H(t,u)}{u}=b(t,u)+[1-q_{+}(t)]\cdot\frac{h(t)}{t},\qquad t\geq t_0,
\end{eqnarray*}
for a continuous function $b:[t_0,+\infty)\times\mathbb{R}\rightarrow\mathbb{R}$ subjected to certain integral restrictions, see the cited
literature, it is desirable to look for positive solutions of the $q_{\pm}$ -- problem (\ref{main_bvp}) with the functions $q_{\pm}$ described
by
\begin{eqnarray}
q_{\pm}(t)=o(1)\qquad\mbox{when }t\rightarrow+\infty.\label{decay_q}
\end{eqnarray}
This will allow for a "free of restrictions" part of $\frac{H(t,u)}{u}$ as large as possible by such an approach.

In this note, using a special result about (\ref{main_bvp}), a flexible criterion for the existence of positive solutions to (\ref{main_PDE})
that decay to $0$ as $\vert x\vert\rightarrow+\infty$ is established. It extends and complements the conclusions of \cite{DMM,EM}.

\section{A boundary value problem for the logarithmic derivative: statement and application}
Let us consider the problem
\begin{eqnarray}
\left\{
\begin{array}{ll}
u^{\prime\prime}+F(t,u)=0,\qquad t\geq t_0>0,\\
\alpha(t)\leq\frac{u^{\prime}(t)}{u(t)}\leq\beta(t),\qquad t\geq t_0,\\
\int_{t_1}^{t_2}\frac{F(s,u(s))}{u(s)}ds\leq\gamma(t_1,t_2),\qquad (t_1,t_2)\in\Gamma,\\
u(t_0)=u_0>0,
\end{array}
\right.\label{bvp}
\end{eqnarray}
where the functions $F:[t_0,+\infty)\times\mathbb{R}\rightarrow[0,+\infty)$, $\alpha$, $\beta:[t_0,+\infty)\rightarrow[0,+\infty)$ and
$\gamma:\Gamma=\{(t_1,t_2):t_1\geq t_0,\thinspace t_{2}\in[t_1,t_1+p]\}\rightarrow[0,+\infty)$ are continuous. Here, $p>0$ is fixed.

It is assumed that
\begin{eqnarray}
\lim\limits_{t\rightarrow+\infty}\alpha(t)=\lim\limits_{t\rightarrow+\infty}\beta(t)=0,\qquad\alpha,\beta\in
L^{2}((t_0,+\infty),\mathbb{R}).\label{master_restr_alphabeta}
\end{eqnarray}
Also, for every $\varepsilon>0$ there exists $\zeta(\varepsilon)\in(0,p)$ such that
\begin{eqnarray}
\gamma(t_1,t_2)<\varepsilon,\qquad 0\leq t_2-t_1\leq\zeta(\varepsilon),\thinspace (t_1,t_2)\in\Gamma.\label{restr_alphabeta2}
\end{eqnarray}

The particular case of interest herein is given by
\begin{eqnarray}
\alpha(t)=\frac{q_{-}(t)}{t}\qquad\mbox{and}\qquad\beta(t)=\frac{q_{+}(t)}{t},\label{partic_alphabeta}
\end{eqnarray}
together with
\begin{eqnarray}
F(t,u)=\left\vert\frac{H(t,u)}{u}-[1-q_{+}(t)]\cdot\frac{h(t)}{t}\right\vert u\label{partic_F}
\end{eqnarray}
for all $t\geq t_0$ and $u\in\cal{C}$ (to be defined later in this section).

We notice that, given $\lambda\in(0,1)$, condition (\ref{decay_q}) yields the existence of $t_{\lambda}\geq t_0$ such that $\sup\limits_{s\geq
t_{\lambda}}q_{+}(s)\leq\lambda$. Consequently, we have
\begin{eqnarray}
\exp\left(\int_{t_0}^{t}\frac{q_{+}(s)}{s}ds\right)=O(t^{\lambda})=o(t)\qquad\mbox{when }t\rightarrow+\infty.\label{decay_u}
\end{eqnarray}
In particular, the solution $u$ of (\ref{bvp}) for $\alpha$, $\beta$ given by (\ref{partic_alphabeta}) will obey the long-time law
\begin{eqnarray*}
\lim\limits_{t\rightarrow+\infty}\frac{u(t)}{t}=0,
\end{eqnarray*}
thus being a solution of (\ref{main_bvp}).

The boundary value problem (\ref{bvp}) is, bluntly speaking, about the existence of a positive solution to a nonlinear ordinary differential
equation that has prescribed long-time behavior for its logarithmic derivative. In this way, the present problem is in the spirit of the
investigation from \cite{M2005}.

The integral condition in the statement of (\ref{bvp}) has a very simple particular case illuminating its presence:
\begin{eqnarray*}
0<\frac{F(t,u)}{u}\leq k<+\infty\quad\mbox{when } u\geq u_{0}\exp\left(\int_{t_0}^{t}\alpha(s)ds\right),\thinspace t\geq t_0,
\end{eqnarray*}
where $k$ is a constant. Here, $\gamma(t_1,t_2)=k(t_{2}-t_{1})$ for all $(t_1,t_2)\in\Gamma$. Since the problem (\ref{bvp}) describes a
non-oscillatory solution to the equation
\begin{eqnarray*}
u^{\prime\prime}+F(t,u)=0,\qquad t\geq t_0,
\end{eqnarray*}
while the comparison equation
\begin{eqnarray*}
u^{\prime\prime}+ku=0,\qquad t\geq t_0,
\end{eqnarray*}
is oscillatory, the analysis in \cite[Theorem 5]{Wong} establishes that the presence in (\ref{bvp}) of the condition involving $\gamma$ does not
diminish in any way the applicability of our technique to the problem presented in the introduction.

The main hypotheses regarding $F(t,u)$ are given following the lines of the Hale-Onuchic integration theory \cite{HaleOnuchic}. By introducing
the set  $\cal{B}$ of functions $b\in C([t_0,+\infty),\mathbb{R})$ which obey for all $t\geq t_0$ the inequalities
\begin{eqnarray}
\alpha(t)\leq b(t)\leq\beta(t),\label{BbB}
\end{eqnarray}
we ask that $F$ be confined by
\begin{eqnarray}
B_{-}(t)&\leq&\frac{1}{u_0}\int_{t}^{+\infty}F\left(s,u_{0}\exp\left(\int_{t_0}^{s}b(\tau)d\tau\right)\right)\nonumber\\
&\times& \exp\left(-\int_{t_0}^{s}b(\tau)d\tau\right)ds\leq B_{+}(t),\qquad t\geq t_0,\label{restr_haleonuchic}
\end{eqnarray}
where
\begin{eqnarray*}
B_{-}(t)=\alpha(t)-\int_{t}^{+\infty}[\alpha(s)]^{2}ds,\qquad B_{+}(t)=\beta(t)-\int_{t}^{+\infty}[\beta(s)]^{2}ds,
\end{eqnarray*}
and also by
\begin{eqnarray}
&&\frac{1}{u_0}\int_{t_1}^{t_2}F\left(s,u_{0}\exp\left(\int_{t_0}^{s}b(\tau)d\tau\right)\right)\nonumber\\
&&\times\exp\left(-\int_{t_0}^{s}b(\tau)d\tau\right)ds\leq\gamma(t_1,t_2),\qquad (t_1,t_2)\in\Gamma,\label{F_equicont}
\end{eqnarray}
for all $b\in\cal{B}$.

In the particular cases of $F(t,u)=A(t)u$ and $F(t,u)=A(t)u^{\sigma}$, where $\sigma\in(0,1)$, the restriction (\ref{restr_haleonuchic}) reads
as
\begin{eqnarray}
B_{-}(t)\leq\int_{t}^{+\infty}A(s)ds\leq B_{+}(t)\label{restr_lin}
\end{eqnarray}
and respectively as
\begin{eqnarray}
B_{-}(t)&\leq& u_{0}^{\sigma-1}\int_{t}^{+\infty}A(s)\nonumber\\
&\times&\exp\left(-(1-\sigma)\int_{t_0}^{s}b(\tau)d\tau\right)ds\leq B_{+}(t)\label{restr_emdenFowler}
\end{eqnarray}
throughout $[t_0,+\infty)$. An immediate simplification of (\ref{restr_emdenFowler}), by means of (\ref{BbB}), is given via the system of
inequalities
\begin{eqnarray*}
\left\{
\begin{array}{ll}
B_{-}(t)\leq u_{0}^{\sigma-1}\int_{t}^{+\infty}A(s)\exp\left(-(1-\sigma)\int_{t_0}^{s}\beta(\tau)d\tau\right)ds\\
B_{+}(t)\geq u_{0}^{\sigma-1}\int_{t}^{+\infty}A(s)\exp\left(-(1-\sigma)\int_{t_0}^{s}\alpha(\tau)d\tau\right)ds.
\end{array}
\right.
\end{eqnarray*}
Here, in the first particular case of $F$, we have
\begin{eqnarray*}
\gamma(t_1,t_2)=\int_{t_1}^{t_2}A(s)ds,
\end{eqnarray*}
while in the second particular case
\begin{eqnarray*}
\gamma(t_1,t_2)=u_{0}^{\sigma-1}\int_{t_1}^{t_2}A(s)\exp\left(-(1-\sigma)\int_{t_0}^{s}\alpha(\tau)d\tau\right)ds
\end{eqnarray*}
for all $(t_1,t_2)\in\Gamma$.

An existence result for problem (\ref{bvp}) now reads as follows.

\begin{theorem}\label{th_ODE} Suppose that the functions $F$, $\alpha$, $\beta$, $\gamma$ satisfy the conditions (\ref{master_restr_alphabeta}), (\ref{restr_alphabeta2}), (\ref{restr_haleonuchic}) and (\ref{F_equicont}). Then, the problem (\ref{bvp}) has a solution $u$. In particular, its logarithmic
derivative lies in $\cal{B}$.
\end{theorem}

Its proof is presented in section \ref{prf}.

To apply the conclusions of Theorem \ref{th_ODE} to the analysis of equation (\ref{main_PDE}), let us introduce the functions $q_{-}(t)\equiv0$ and
$q_{+}$ subjected to (\ref{decay_q}) and (\ref{decay_u}). Here, $t_0=(n-2)R^{n-2}$, $R>0$. Define also the functions $H$, $h$ on the basis of
(\ref{def_H}) and assume, for simplicity, that
\begin{eqnarray}
\frac{H(t,u(t))}{u(t)}\geq[1-q_{+}(t)]\cdot\frac{h(t)}{t},\qquad u\in{\cal{C}},\label{comp_Ha}
\end{eqnarray}
throughout $[t_0,+\infty)$. The set $\cal{C}$ consists of all the functions $u\in C^1([t_0,+\infty),(0,+\infty))$ with $u(t_0)=u_0$ and
$\frac{u^{\prime}}{u}\in\cal{B}$. We fix $u_0>0$ such that $\frac{u_0}{t_0}\leq\varsigma$.

The condition (\ref{comp_Ha}) says, practically, that $m$ from (\ref{core}) is in a certain sense larger than $g$. We recall that \cite[Theorem
1, Remark 1]{EM} dealt with the complementary case, namely $a(r)\leq lr^{n-2}g(r)$ for all $r\geq R$, where $m(r,U)=a(r)w(U)$ and $l>0$ is a
constant.

Our main contribution here is the next result.

\begin{theorem}\label{th_PDE} Assume that $g(r)\geq 0$ for all $r\geq R$. Suppose further that, when $r\geq R$ and $R\leq r_{1}\leq r_{2}\leq
\left(r_{1}^{n-2}+\frac{p}{n-2}\right)^{\frac{1}{n-2}}$, one has
\begin{eqnarray*}
\int_{r}^{+\infty}M(\tau,u((n-2)\tau^{n-2}))d\tau\leq B_{+}((n-2)r^{n-2})
\end{eqnarray*}
and
\begin{eqnarray*}
\int_{r_1}^{r_2}M(\tau,u((n-2)\tau^{n-2}))d\tau\leq \gamma((n-2)r_{1}^{n-2},(n-2)r_{2}^{n-2})
\end{eqnarray*}
for all $u\in\cal{C}$, where
\begin{eqnarray*}
M(\tau,u)=\frac{\tau}{n-2}\left\{\frac{m\left(\tau,\frac{u}{(n-2)\tau^{n-2}}\right)}{u}-[1-q_{+}((n-2)\tau^{n-2})]\frac{g(\tau)}{\tau^{n-2}}\right\}.
\end{eqnarray*}

Then, the equation (\ref{main_PDE}) has a positive solution $u(x)$, defined in $G_{R}$, such that $\lim\limits_{\vert
x\vert\rightarrow+\infty}u(x)=0$.
\end{theorem}

Its proof is presented in the next section.

\section{Proofs}\label{prf}

\textit{Proof of Theorem \ref{th_ODE}}. The existence of a solution to the problem (\ref{bvp}) will be demonstrated by employing the Schauder fixed point theorem.

Let $A(t_{0})$ be the real linear space of continuous functions $b:[t_{0},+\infty)\rightarrow\mathbb{R}$ which satisfy
\begin{eqnarray*}
\lim_{t\rightarrow +\infty}b(t)=l_{b}\in \mathbb{R}.
\end{eqnarray*}
If endowed with the standard sup-norm $\Vert\cdot\Vert$, $A(t_{0})$ becomes a Banach space. Avramescu's criterion \cite{A} for relative
compactness of subsets in this space asks from $S\subset A(t_{0})$ to be norm-bounded, \textit{equicontinuous} (meaning that all the functions
from $S$ are uniformly continuous in the same way) and \textit{equiconvergent} (that is, all the functions $b\in S$ approach their limits
$l_{b}$ in an uniform way) in order to be relatively compact.

We start by noticing that ${\cal{B}}\subset A(t_0)$ with $l_{b}=0$ for all $b\in\cal{B}$. Since $\Vert b\Vert\leq\Vert\beta\Vert$ throughout
$\cal{B}$, the set is bounded. It is also easy to conclude that it is convex and closed in the norm topology of $A(t_0)$.

Further, we define the integral operator $T:{\cal{B}}\rightarrow A(t_0)$ through the formula
\begin{eqnarray*}
T(b)(t)&=&\int_{t}^{+\infty}[b(s)]^{2}ds+\frac{1}{u_0}\int_{t}^{+\infty}F\left(s,u_{0}\exp\left(\int_{t_0}^{s}b(\tau)d\tau\right)\right)\\
&\times&\exp\left(-\int_{t_0}^{s}b(\tau)d\tau\right)ds,\qquad t\geq t_0.
\end{eqnarray*}

The assumptions (\ref{BbB}), (\ref{restr_haleonuchic}) imply that
\begin{eqnarray}
T({\cal{B}})\subseteq\cal{B}.\label{welldef_T}
\end{eqnarray}

We shall establish that $T:{\cal{B}}\rightarrow{\cal{B}}$ is continuous.

Set $\varepsilon >0$. There exists $T_{\varepsilon}\geq t_0$ such that $\beta(T_{\varepsilon})\leq\frac{\varepsilon}{6}$. Introduce also the
functions
\begin{eqnarray*}
g(t,w)=\left(u_{0}e^{w}\right)^{-1}F\left(t,u_{0}e^{w}\right)\qquad\mbox{and}\qquad x(t;b)=\int_{t_{0}}^{t}b(\tau)d\tau
\end{eqnarray*}
for all $w\geq0$ and $b\in{\cal{B}}$.

Since $g:[t_{0},T_{\varepsilon}]\times[0,w_{\varepsilon}]\rightarrow {\R}$ is uniformly continuous, where
$w_{\varepsilon}=u_0\exp\left(\int_{t_0}^{T_{\varepsilon}}\beta(\tau)d\tau\right)$, there exists $\delta (\varepsilon)>0$ such that
\begin{eqnarray*}
\vert g(s,w_{1})-g(s,w_{2})\vert\leq\frac{\varepsilon}{3T_{\varepsilon}},\qquad s\in [t_{0},T_{\varepsilon}],
\end{eqnarray*}
for all $w_{1,2}\in[0,w_{\varepsilon}]$ with $\vert w_{1}-w_{2}\vert\leq\delta (\varepsilon)$. Notice that
\begin{eqnarray*}
\vert x(t;b_{1})-x(t;b_{2})\vert\leq\int_{t_{0}}^{t}\vert b_{1}(\tau)-b_{2}(\tau)\vert d\tau\leq T_{\varepsilon}\Vert b_{1}-b_{2}\Vert.
\end{eqnarray*}

Now, given $b_{1,2}\in\cal{B}$ with $\Vert
b_1-b_2\Vert\leq\eta(\varepsilon)=\min\left\{\frac{\varepsilon}{6T_{\varepsilon}\Vert\beta\Vert},\frac{\delta(\varepsilon)}{T_{\varepsilon}}\right\}$,
we have the estimates
\begin{eqnarray*}
&&\vert T(b_1)(t)-T(b_2)(t)\vert\\
&&\leq\int_{t_0}^{T_{\varepsilon}}\vert [b_{1}(s)]^{2}-[b_{2}(s)]^{2}\vert
ds+2\int_{T_{\varepsilon}}^{+\infty}[\beta(s)]^{2}ds\\
&&+\int_{t_0}^{T_{\varepsilon}}\vert g(s,x(s;b_1))-g(s,x(s;b_2))\vert ds+2B_{+}(T_{\varepsilon})\\
&&\leq\int_{t_0}^{T_{\varepsilon}}[b_{1}(s)+b_{2}(s)]ds\cdot\Vert
b_{1}-b_{2}\Vert+2\int_{T_{\varepsilon}}^{+\infty}[\beta(s)]^{2}ds\\
&&+\frac{\varepsilon}{3}+2B_{+}(T_{\varepsilon})\\
&&\leq2\Vert\beta\Vert\cdot T_{\varepsilon}\Vert b_{1}-b_{2}\Vert+2\beta(T_{\varepsilon})+\frac{\varepsilon}{3}\leq\varepsilon,
\end{eqnarray*}
that is, $\Vert T(b_1)-T(b_2)\Vert\leq\varepsilon$.

In order to apply Schauder's fixed point theorem to operator $T$, it remains to prove that the set $T({\cal{B}})$ verifies the hypotheses of
Avramescu's criterion.

The first one, namely the boundedness of $T({\cal{B}})$, follows from (\ref{welldef_T}).

The equicontinuity of $T({\cal{B}})$ is a consequence of the estimate
\begin{eqnarray*}
0\leq T(b)(t_1)-T(b)(t_2)\leq\Vert\beta\Vert^{2}(t_2-t_1)+\gamma(t_1,t_2),\qquad(t_1,t_2)\in\Gamma,
\end{eqnarray*}
together with (\ref{restr_alphabeta2}).

The equiconvergence of $T({\cal{B}})$ is implied by
\begin{eqnarray*}
0\leq T(b)(t)=\vert T(b)(t)-l_{T(b)}\vert\leq\beta(t)=o(1)\qquad\mbox{as }t\rightarrow+\infty
\end{eqnarray*}
for all $b\in\cal{B}$.

The solution of problem (\ref{bvp}) has the formula
\begin{eqnarray*}
u(t)=u_0\exp\left(\int_{t_0}^{t}b_{0}(s)ds\right),\qquad t\geq t_0,
\end{eqnarray*}
where $b_0\in\cal{B}$ is the fixed point of operator $T$. $\square$

The following lemma is of use for Theorem \ref{th_PDE}.

\begin{lemma}\label{lemma:subsup}
\emph{ (see \cite{NS})} If there exist a nonnegative subsolution $w$ and a positive supersolution $v$ to (\ref{main_PDE}) in $G_R$, such that
$w(x)\leq v(x)$ for $x\in \overline{G}_R$, then (\ref{main_PDE}) has a solution $u$ in $G_R$ such that $w\leq u\leq v$ throughout
$\overline{G}_R$. In particular, $u=v$ on $\vert x\vert=R$.
\end{lemma}

\textit{Proof of Theorem \ref{th_PDE}}. Consider the positive, twice continuously differentiable functions given by
\begin{eqnarray*}
U(x)=y(r)=\frac{u_{1}(t)}t,\qquad t\geq t_0,
\end{eqnarray*}
where $r=\vert x\vert=\theta (t)$. Here, $u_{1}$ is the solution of problem (\ref{bvp}) obtained at Theorem \ref{th_ODE}. Since the range of
$q_{+}$ is a subset of $[0,1]$, we have
\begin{eqnarray*}
\frac{u(t)}{t}\leq\frac{u_0}{t}\exp\left(\int_{t_0}^{t}\frac{q_{+}(s)}{s}ds\right)\leq\frac{u_0}{t_0}\leq\varsigma,\qquad t\geq t_0,
\end{eqnarray*}
for all $u\in\cal{C}$. This estimate allows us to use the comparison inequality (\ref{core}) in the following.

By a straightforward computation we get that
\begin{equation}
t\theta^{\prime}(t)=\frac1{n-2}\theta (t)\label{set_1}
\end{equation}
and
\begin{equation}
\left\{
\begin{array}{ll}
\frac{du_1}{dt}=y+t\theta^{\prime}(t)\frac{dy}{dr} \\
\frac{d^2u_1}{dt^2}=\frac{n-1}{n-2}\theta^{\prime}(t) \frac{dy}{dr}+\frac{\theta(t)\theta^{\prime}(t)} {n-2}\frac{d^2y}{dr^2}\cdot
\end{array}
\right.\label{set_2}
\end{equation}

Further, taking into account (\ref{set_1}) and (\ref{set_2}), we have
\begin{eqnarray*}
&&r^{n-1}\left(\Delta U+f(x,U)+g(\vert x\vert )x\cdot \nabla U\right)\\
&&=\frac d{dr}\left(r^{n-1}\frac{dy}{dr}\right)
+r^{n-1}f(x,U)+r^ng(r)\frac{dy}{dr}\nonumber\\
&&=\frac{n-2}{\theta (t)\theta^{\prime}(t)}[\theta (t)]^{n-1}\left[
u_{1}^{\prime\prime}(t)+\frac 1{n-2}\theta (t)\theta^{\prime}(t)f(x,U)\right.\\
&&+\left.\theta (t)\theta^{\prime}(t)g(\theta(t))\left(u_{1}^{\prime}(t) -\frac{u_1(t)}t\right)\right],
\end{eqnarray*}
for any $t\geq t_0$.

We have obtained that
\begin{eqnarray*}
&&\vert x\vert^{n-1}\left(\Delta U+f(x,U)+g(\left|x\right|)x\cdot\nabla
U\right)\\
&&\leq \frac{n-2}{\theta (t)\theta^{\prime}(t)}[\theta(t)]^{n-1}\left[ u_{1}^{\prime\prime}(t)+F(t,u_{1}(t))\right]=0,
\end{eqnarray*}
where the function $F$ is given by (\ref{partic_F}).

Now, $U$ is a positive super-solution of (\ref{main_PDE}). Also, the trivial solution of (\ref{main_PDE}) is its (nonnegative) sub-solution.
According to Lemma \ref{lemma:subsup}, there exists a nonnegative solution $u$ to (\ref{main_PDE}), defined in $\overline{G}_{R}$. Since
\begin{eqnarray*}
(\Delta+g(\vert x\vert)x\cdot\nabla)(-u)=f(x,u)\geq 0,
\end{eqnarray*}
the strong maximum principle (\cite{F00}) can be applied to $-u$. This means that the function $-u$ cannot attain a nonnegative maximum at a
point of $G_R$ unless it is constant. Since $-u$ is negative on $\{x:\vert x\vert=R\}$ and $-u(x)\leq 0$ throughout $\overline{G}_R$ as $u$ is
confined between $0$ and a positive super-solution $U$, it follows that $-u$ cannot have zeros.

We conclude that $u$ is a positive solution of (\ref{main_PDE}) that decays to $0$ when $\vert x\vert\rightarrow+\infty$. Furthermore, via
(\ref{decay_u}), we can compute the decay rate of $u$:
\begin{eqnarray*}
0<u(x)\leq\frac{u_{1}((n-2)\vert x\vert^{n-2})}{(n-2)\vert x\vert^{n-2}}=O\left(\vert x\vert^{(\lambda-1)(n-2)}\right)\quad\mbox{when }\vert
x\vert\rightarrow+\infty.
\end{eqnarray*}

The proof is complete. $\square$

\section{Conclusion}

To emphasize the significance of Theorem \ref{th_PDE}, let us consider a particular case of the comparison function $m$ from (\ref{core}), namely
\begin{eqnarray}
m(\vert x\vert,U)=a(\vert x\vert)U,\qquad x\in G_R,\thinspace U\in[0,\varsigma],\label{part_m}
\end{eqnarray}
for a continuous function $a:[R,+\infty)\rightarrow[0,+\infty)$.

In some of the recent literature \cite{Constantin1996,Constantin1997,Constantin2005,E,Yin} regarding (\ref{main_PDE}), the leading hypothesis was
\begin{eqnarray}
\int_R^{+\infty}r\left[a(r)+\left|g(r)\right|\right]dr<+\infty.\label{class_ref}
\end{eqnarray}
The conclusion of these papers reads, practically, as follows: since the
elliptic equation
\begin{equation}
\Delta u+f(x,u)=0,\qquad \left|x\right|>R,\label{unperturbed}
\end{equation}
has a positive solution decaying to $0$ as $\vert x\vert\rightarrow+\infty$ under the hypotheses
(\ref{core}), (\ref{part_m}) and $\int_{R}^{+\infty}ra(r)dr<+\infty$, its "small" perturbation by the term
"$g(\left| x\right|)x\cdot \nabla u$", namely (\ref{main_PDE}),
where the degree of "smallness" is given by (\ref{class_ref}), will
preserve this feature.

The hypotheses of Theorem \ref{th_PDE}, if verified,
reveal that the behavior of certain solutions to (\ref{main_PDE})
for large $\left|x\right|$'s instead of being controlled by this
summing action of functions $a$, $g$ is actually governed by their
coupling. In our case, that is, it might happen that
for certain functions $a$, $g$, where
\begin{eqnarray*}
0\leq g(r)\leq \frac{a(r)}{n-2},\qquad r\geq R,
\end{eqnarray*}
and
\begin{eqnarray*}
\int_{r}^{+\infty}n(\tau)d\tau\leq B_{+}((n-2)r^{n-2}),\qquad r\geq R,
\end{eqnarray*}
and
\begin{eqnarray*}
\int_{r_1}^{r_2}n(\tau)d\tau\leq\gamma((n-2)r_{1}^{n-2},(n-2)r_{2}^{n-2}),
\end{eqnarray*}
where $R\leq r_{1}\leq r_{2}\leq
\left(r_{1}^{n-2}+\frac{p}{n-2}\right)^{\frac{1}{n-2}}$ and
\begin{eqnarray*}
n(r)=\frac{r^{3-n}}{n-2}\left\{\frac{a(r)}{n-2}-[1-q_{+}((n-2)r^{n-2})]g(r)\right\},\qquad r\geq R,
\end{eqnarray*}
the unperturbed equation
(\ref{unperturbed}) does not have any vanishing at $+\infty$
solution besides the trivial solution.

\textbf{Acknowledgment} We are grateful to a referee for his/hers useful suggestions.

\end{document}